\documentclass[a4paper,11pt]{article}
\usepackage{pstricks}
\usepackage{amsmath}
\usepackage{amsfonts}
\usepackage{amscd}
\usepackage{amssymb}
\usepackage[all]{xy}

\newcommand{\al}{{\alpha}}

\newcommand{\Om}{{\Omega}}
\newcommand{\om}{{\omega}}
\newcommand{\eps}{{\varepsilon}}

\newcommand{\la}{{\lambda}}

\newcommand{\mapto}[1]{\buildrel#1\over \longrightarrow}

\newcommand{\rar}{\rightarrow}

\newcommand{\x}{\times}
\newcommand{\er}{\mathbb R}

\newcommand{\hook}{\hookrightarrow}

\newcommand{\ifff}{if and only if\ }

\newcommand{\NI}{{\noindent}}

\newcommand{\bk}{\hfill $\Box$}

\newtheorem{theorem}{Theorem}[section]
\newtheorem{thm}[theorem]{Theorem}
\newtheorem{cor}[theorem]{Corollary}
\newtheorem{defin}[theorem]{Definition}
\newtheorem{rem}[theorem]{Remark}
\newtheorem{lemma}[theorem]{Lemma}
\newtheorem{prop}[theorem]{Proposition}
\newtheorem{ex}[theorem]{Example}
\newtheorem{que}[theorem]{Question}

\numberwithin{equation}{section}

\title{\bf Presymplectic manifolds}
\author{Bogus\l\/aw Hajduk $^*$ and Rafa\l\ Walczak \thanks{Both authors
were partially supported by grant of the Ministry of Science and
Higher Education of Poland, no. 1P03A 03330} }

\date{}

\begin{document}

\maketitle

\begin{abstract} A presymplectic structure on an odd-dimensional
manifold is given by a closed 2-form which is non-degenerate, i.e.,
of maximal rank. We investigate geometry of presymplectic manifolds.
Some basic theorems analogous to those in symplectic and contact
topology are given and applied to study constructions of
presymplectic manifolds. In particular, we describe how to glue
presymplectic manifolds along isomorphic presymplectic submanifolds,
including surgery on presymplectic circles.

\noindent {\bf Keywords:} presymplectic manifold,  differential form
of maximal rank, contact form

\noindent {\bf AMS classification(2000)}: Primary 53D05; Secondary
57S25.

\noindent {\bf AMS classification(1991)}: Primary 53C15; Secondary
57S25.

\end{abstract}

\bigskip

\section{Introduction}

Investigation  of symplectic and contact structures is an important
and developing part of geometry of manifolds. In this paper we study
presymplectic manifolds which mediate between those two
types of structures.

A {\it presymplectic form} is a closed differential 2-form  of
maximal rank on an odd-dimensional manifold. Thus if the manifold
$M$ is of dimension  $2n+1,$ $\om$ is a presymplectic form on $M,$
then the rank of the form is $2n.$ This means that at each point
$x\in M$ there exists a (non-unique!) codimension one subspace of
$T_xM$ such the form is a symplectic linear form on it. This gives a
linear symplectic, hence admitting a complex structure, subbundle
and a unique subbundle $\mathcal R=\{ V\in T_xM| \iota_V\om = 0\}$
of dimension one, which we call {\it Reeb bundle.} Throughout the
paper we assume that $M$ is an oriented manifold. Then there is an
orientation of $\mathcal R$ given by the orientation of $M$ and the
orientation of the symplectic subbundle provided by $\om .$ In
particular, $TM$ has a reduction of the structure group to $U(n)\x
\mathbf{1}.$ Furthermore, the Reeb bundle $\mathcal R$ defines a
foliation of dimension one and, for any choice of non-zero section
$R$ of $\mathcal R,$ a flow without fixed points.
 We shall call
the  foliation the {\it Reeb foliation} of the presymplectic form
$\om ,$ and any such $R$  a {\it Reeb vector field.}

There is a natural source of presymplectic structures.

\begin{lemma}\label{cod1} If $i:Q\subset M$ is a codimension one submanifold of a
symplectic manifold $(M,\Om ),$ then $i^*\Om$ is a presymplectic
form on $Q.$
\end{lemma}

As we noticed above, if an oriented manifold $M$ admits a
presymplectic form, then the tangent bundle of $M$ has a reduction
to $U(n)\x \bold 1\subset O(2n+1)$ defined by Reeb subbundle and a
choice of a complex structure in a complementary  symplectic
subbundle. As in the case of almost complex structures defined by a
symplectic form, there is a large family of possible choices. Any
such choice is provided by a bundle map $J:TM\rar TM$ such that
$J(\mathcal R) = 0$ and $J^2=-Id$ on the subbundle $Im\, J.$ One can
see that existence of such a reduction is equivalent to existence of
a 2-form of maximal rank (not necessarily closed), cf. \cite{B}.
This is analogous to existence of an almost complex structure for a
symplectic form.

However,  existence  of a 2-form of maximal rank is
also sufficient for the existence of presymplectic form. It reduces
the problem of existence to homotopy theoretical questions due to the following fundamental
theorem.

\begin{thm}\label{pretheorem} If $\mathcal{S}_{non-deg}(M)$ is the space of all
2--forms of maximal rank on a closed manifold $M^{2n+1}$ and
$\mathbb{S}_{presymp}(M;a)$ is the space of all presymplectic forms
in a given cohomology class $a,$ then
$$\mathbb{S}_{presymp}(M,a) \hook \mathcal{S}_{non-deg}(M)$$ is a homotopy
equivalence.
\end{thm}

Main technical novelty of the present note is that the assumption on
the cohomology class can be removed and a relative case also holds,
cf. Section \ref{symplhom}. This is later applied to various
constructions.

The first proof that there is a bijective correspondence of
connected components, with a relative version,  was given by Dusa
McDuff in \cite{MD}.

\begin{thm}\label{dusa} Any 2-form $\om$ of maximal rank can be
deformed in the space of forms of maximal rank, to a presymplectic
form $\om'.$ If $\om$ is presymplectic in a neighborhood of a
compact set $K,$ then there exists a deformation which is constant
in a neighborhood of $K.$
\end{thm}

 With some
effort this can be extended to prove the homotopy equivalence. Later
Eliashberg and Mishachev \cite{EM} gave a proof based on Gromov's
existence theorem of symplectic forms on open manifolds. The idea
can be easily explained:  $\om$ gives a form $\om + \rho \wedge dt,\
\rho (R)=1,$ of maximal rank on $M\x \mathbb R$ (not closed in
general), by Gromov's theorem it can be deformed to a symplectic
form. Its restriction to $M\x \{ 0\}$ yields a presymplectic form by
Lemma \ref{cod1}.

As an immediate corollary we get that any parallelizable manifold
admits a presymplectic form, e.g. any oriented 3-manifold is
presymplectic.

The fact that it is relatively simple to construct a presymplectic
forms can be useful for contact and symplectic topology. For
example, one can ask the following.

\begin{que} \label{econtact} Does any closed presymplectic
manifold admit a contact form?
\end{que}

A rather natural idea is to consider a presymplectic form of zero
cohomology class, $\om = d\al ,$ and look for a deformation of $\al$
to a contact form. By \cite{G}, this can be done in dimension 5 for
any closed simply connected presymplectic manifold. However, except
for this result, the classical case of dimension 3 and a number of
special cases, no answers to this question are known.  In section
\ref{contacttype} we show that vanishing of basic cohomology class
with respect to its Reeb foliation is an obstruction for  a
presymplectic form to be of {\it contact type} (i.e. equal to
differential of a contact form).

Consider two other problems, which proved to be very difficult and
 both are open in dimensions greater than four. First,
let $M$ be an odd-dimensional closed manifold.

\begin{que} \label{taubes} Is it true that if $M\x S^1$ is symplectic, then
$M$ fibres over $S^1?$ \end{que}

For $M$ of dimension three this question was posed by Taubes and
answered positively, after a series of partial results of many
authors, by Friedl and Vidussi \cite{FV}. Their proof uses Seiberg
-Witten invariants and it does not extend to higher dimensions.

 A closed  manifold $X$ of dimension $2n$
 we shall call {\it homotopically symplectic} if $X$ is almost complex and there exists
 $x\in H^2(X,\mathbb R)$ such that
 $x^n\neq 0.$
These two conditions are the only known necessary conditions for existence of a symplectic
 structure for a general closed manifolds of dimension greater than
 four.

Note that the term {\it c-symplectic} or {\it  cohomologically
symplectic} is used for manifolds which satisfy  the cohomological
part of the condition.

\begin{que}\label{esymp} Does any closed, homotopically symplectic manifold
admit a symplectic form?
\end{que}

Questions ref{taubes,esymp} cannot simultaneously have positive
answers. Namely, there are presymplectic manifolds which do not
fibre over the circle, but their products with the circle satisfy
assumption of \ref{esymp}, e.g. connected sum of two copies of tori
$T^{2k+1}\# T^{2k+1}.$

The principal aim of the present paper is to give some basic
theorems on presymplectic manifolds and to provide means to
constructions of presymplectic manifolds. It includes Moser type
theorems, tubular neighborhoods of presymplectic submanifolds and
constructions of presymplectic manifolds. For symplectic and contact
manifolds such theorems open a way to analyze their topology. In the
presymplectic case the Reeb foliation plays an important role. For
example, on a closed symplectic manifold we have Moser's theorem
saying that a path $\om_t$ of symplectic forms with constant
cohomology class is given by an isotopy, i.e., there exists an
isotopy $f_t$ such that $\om_t=f_t^*\om_0,$ can not be extended to
presymplectic manifolds without an additional assumption on the Reeb
foliation. This assumption is, roughly speaking,  that the basic
cohomology class of the presymplectic form is constant. A similar
problems appear if one analyzes behavior of a presymplectic form in
a neighborhood of a presymplectic submanifold. For tubular
neighborhoods there is no isotopy of different tubular neighborhoods
of a given submanifolds, which allows various constructions in
symplectic and contact topology. Here we have only a weaker
statement of a presymplectic isomorphisms of such neighborhoods.

Despite of the fact that existence of presymplectic form boils down
to homotopy theoretical, hence "soft" question, some effective
construction methods can be useful. We describe first a  Thurston
type construction of presymplectic forms on bundles in Section
\ref{thurston}.

  Another construction which may have
interesting application gives a presymplectic structure on any
 presymplectic open book decomposition.
This technique is  well known in contact theory. By results of
Giroux and Mohsen \cite{GM}, any closed contact manifold admits a
compatible open book decomposition and in dimension three this  was
successfully used to some classification problems. For presymplectic
manifolds the construction of open book, based on Donaldson method
of quasi-holomorphic sections, can be repeated, cf. \cite{MMP}. We
work in reverse direction, namely for a given open book
decomposition we give sufficient condition to construct a
presymplectic form with  Reeb bundle transversal to pages.

In fact, we do this in more general case of {\it star-like
structures.} This is a straightforward generalization of the notion
of open book decomposition to the case when the complement of a
submanifold of codimension $2k$ admits a symplectic fibration over
$S^{2k-1}.$

Finally, we describe how to glue two presymplectic manifolds along
tubular neighborhoods of isomorphic  presymplectic submanifolds.  A
particular case of 1-dimensional submanifold, hence a closed leaf of
the Reeb foliation, is in fact the  classical surgery on a circle.
We show that in dimension greater than 3 any presymplectic form is
homotopic to one with a closed orbit in each element of a generating
set of the fundamental group, thus one can always kill fundamental
group using such modifications.

In our terminology we follow McDuff \cite{MD}. The term
presymplectic is used in a wider sense in papers on quantization,
e.g. \cite{CKT,KT}, where  a manifold is called presymplectic if it
is endowed with a closed, not necessarily non-degenerate, 2-form.
What we call presymplectic form was called odd-symplectic form in
\cite{Gi} or simply closed form of maximal rank in \cite{EM}. In
hamiltonian mechanics context an exact presymplectic form is called
hamiltonian.

In the sequel we will consider only {\it closed} smooth manifolds and
smooth mappings.

\section{Homotopically presymplectic forms}\label{hompre}

Let $M$ be a closed manifold of dimension $2n+1.$ We say that $M$ is
{\it homotopically presymplectic} if the tangent bundle $TM$ of $M$
admits a reduction to $U(n)\simeq U(n)\x \mathbf{1} \subset
SO(2n+1).$ The name almost contact manifold has also been used for such
structure (see \cite{G1,MMP}). Our terminology attempts to unify
presymplectic and symplectic cases. If $M$ is closed, then we call
it {\it strongly homotopically presymplectic} if it is homotopically
presymplectic and there exists a class $x\in H^2(M)$ such that $x^n$
is non-zero.

\begin{prop} $M\times S^1$ is homotopically symplectic
if and only if $M$ is strongly homotopically presymplectic.
\end{prop}

{\bf Proof.} Assume that $M\times S^1$ is homotopically symplectic.
If $u\in H^2(M\times S^1)$ satisfies $u^{n+1} \neq 0,$ then its
restriction to $H^2(M\times
*)$ satisfies, by the K\"unneth formula, $x^n\neq 0.$ Moreover, any codimension
one subbundle  of a complex bundle has a reduction to $U(n),$ thus
$M$ is strongly homotopically presymplectic.

 If $x \in H^2(M)$ satisfies $x^n \neq 0,$ then there exists a class $\beta \in H^1(M)$
 such that $x^n\cup \beta$ is non-zero, hence  $(x + \beta \cup [dt])^{n+1}\neq 0.$
Furthermore, if $TM$ admits a reduction to $U(n),$ then we can
easily get an almost complex structure $J$ associated with
decomposition $TM=\xi\oplus \eps^1$ by setting
$$J\left(\frac{\partial}{\partial t}\right) = X, \\ JX = -\frac{\partial}{\partial t},$$
where X is a non-zero vector field on $M$ corresponding to the
trivial subbundle  $\epsilon^1$ and $J$ is a complex structure on
$\xi .$ \hfill $\square$

Consider an odd-dimensional manifold $M$ which is strongly
homotopically presymplectic, but do not fibers over the circle, e.g.
$T^{2n-1}\sharp T^{2n-1}.$ Then $M\x S^1$ is strongly homotopically
symplectic. If such manifold is
symplectic, then it is a counter-example to Question \ref{taubes}.

Any symplectic form on $M\x S^1$ gives a loop of presymplectic
forms.   By Lemma \ref{cod1}, there is a continuous mapping from the
space of all symplectic forms on $M \times S^1$ into the space of
(free) loops of presymplectic forms on $M.$ Question \ref{taubes}
leads to the following: does this mapping induce a surjective map
on $\pi_0?$ Equivalently, one can ask if any loop of presymplectic
forms can be deformed to a loop obtained from a symplectic form on
$M\x S^1?$ We show a counter-example on the 3-torus. In particular,
the loop we give is non-contractible.

It is not difficult to construct a loop of non-degenerate 2-forms on
$T^3$ such that the associated almost complex structure $J$ on $T^4$
has non-vanishing first Chern class. It is obtained from a
non-degenerate 2-form on $T^3$ with non-trivial $c_1$ by
restriction. By Theorem \ref{dusa}, we deform this loop to a loop of
presymplectic forms. This loop is not obtained from a symplectic
form, since for any almost complex structure $J$ tamed by a
symplectic structure on $T^4$ we have $c_1(J)=0,$ as proved by
Donaldson.

\section{ Homotopy of presymplectic forms}\label{symplhom}

In the proof of Theorem \ref{pretheorem} due to Mishachev and
Eliashberg, the first step is to pass from a form $\om$ of maximal
rank on $M^{2n+1}$ to a form of maximal rank on $M\x \mathbb R.$
This is very simple.  If $R$ is a Reeb vector field of $\om ,$ then
there exists a 1-form $\eta$ such that $\eta (R)\equiv 1.$ Let
$p:M\x \er\rar M$ be the projection and $\theta$ be the standard
parameter of $\er ,$ so that $d\theta$ is the standard  (translation
invariant) 1-form on $\er .$ Then $p^*\om + d\theta\wedge\eta $ is a
form of maximal rank (non-degenerate) extending $\om$ to $M\x \er .$
However, it can be non-closed even if $\om$ is. The following
observation yields a useful improvement.

\begin{lemma} Let $C\in \er .$ For $\delta$ small enough the form
$\phi = p^*\om + \delta (d\theta\wedge\eta - \theta\, d\eta)$ on
$M\x [-C,C]$ is non-degenerate. If $\om$ is closed in an open set
$U\subset M,$ then $\phi$ is symplectic in $U\x [-C,C]$ and $\phi
|M\x \{ 0\} =\om.$
\end{lemma} \hfill $\square$

Now, after applying the relative version of Gromov's theorem, one comes to
a relative version of Theorem \ref{pretheorem}.

\begin{thm}\label{pretheoremrel} Let $\mathcal{S}_{non-deg}(M,K;a_0)$ be
the space of all 2--forms of maximal rank on a closed manifold
$M^{2n+1}$ which restricts, in a neighborhood of a compact subset
$K\subset M,$ to a presymplectic form in the cohomology class $a_0.$
Let furthermore $\mathbb{S}_{presymp}(M,K;a)$ be the space of all presymplectic
forms in a given cohomology class $a$ which is equal to $a_0$ under
restriction to a neighborhood of $K.$ Then
$$\mathbb{S}_{presymp}(M,K;a) \hook \mathcal{S}_{non-deg}(M,K;a_0)$$
is a homotopy
equivalence.
\end{thm} \hfill $\square$

\begin{cor}\label{path} If  $\eta_t, \ t\in[0,1]$ is a path of
non-degenerate forms on $M$ connecting some presymplectic forms
$\eta_0, \eta_1,$ such that $[\eta_0]=[\eta_1]$ and  $\eta_t$ is
presymplectic in a neighborhood of $K,$ then $\eta_t$ can be
deformed to a path $\tilde\eta_t$ of presymplectic forms which
satisfies:

\begin{enumerate}
\item $\tilde\eta_0=\eta_0,\ \ \tilde\eta_1=\eta_1,$
\item $\tilde\eta_t$ is equal to $\eta_t$ in a neighborhood of $K$
for any $t\in [0,1].$
\end{enumerate}
\end{cor}

In particular, we see that if $M^{2n+1}$ is presymplectic, then for
any two cohomologous presymplectic forms $\om_1,\om_2 $ which can be
connected by a path of non-degenerate forms $\om_{t},$ there is a
path $\om_t'$ of presymplectic forms joining $\om_1$ and $\om_2,$
and $\om_t'$ can be prescribed in a neighborhood of $K.$

We will show now that the assumption that the forms are cohomologous
can be removed. We give a relative version as well.

\begin{thm}\label{prehomotopy} Assume that $\om,\om'$  are two
presymplectic forms on a closed manifold $M^{2n+1}$ and that they are
joined by a path $\{\om_t\}_{t \in [0;1]}$ of non-de\-ge\-ne\-ra\-te
2-forms. Then $\om,\om'$ are homotopic through presymplectic forms.
Furthermore, the following relative version of this theorem also
holds: if $i:K \subset M$ is a compact subset, $\{\om_t\}$ is
presymplectic in an open neighborhood of $K,$ $[i^*\om_t]$ is
constant, then there exists a presymplectic homotopy equal to $\om_t$ in a
neighborhood of $K,$ up to a presymplectic deformation with fixed
ends.
\end{thm}

{\bf Proof}.  Let us define a set $\Lambda_{\om,K}$ as the set of
elements $u\in H^2(M;\er )$ which satisfy the following condition:

(*)  if $[\om']=u$ and there is a path of non-degenerate forms
joining  $\om \ \text{and} \ \om',$ satisfying assumptions of the
relative part of  Theorem \ref{prehomotopy} rel K, then there is a
path of presymplectic forms joining $\om \ \text{and} \ \om'$
equal to the given path in a neighborhood of $K.$

We will show  first that if the condition (*) is satisfied by a
form, then it is satisfied by all other forms in the same cohomology
class. In fact, assume that we have two presymplectic forms
$\om_1,\om_2$ such that $u=[\om_1]=[\om_2]$ and there are  two paths
of non-degenerate forms $\om^1_t,\om^2_t$ joining $\om$ with,
respectively, $\om_1, \om_2$. By Corollary \ref{path}, the path
$\om^1_{-t}\star \om^2_t$ can be  homotoped (rel a neighborhood
$U_K$ of $K$) to a path of presymplectic forms in the class $u$
joining $\om_1$ and $\om_2.$ By assumptions, there is a (rel $U_K$)
path of presymplectic forms joining $\om$ and $\om_1.$ This path in
a neighborhood of $K$ is equal to $\om^1_t\star \om^1_{-t}\star
\om^2_t,$ which is deformable to $\om^2_t.$

Our assumption that the class of $i^*\om_t$ is constant implies that
$\Lambda_{\om,K}\subset [\om]+j^*H^2(M,K,\mathbb{R}),$ where
$j:(M,\emptyset) \hookrightarrow (M,K)$ is the inclusion.

We are going to prove that $\Lambda_{\om,K}$ is both open and closed,
thus equal to   $[\om]+j^*H^2(M,K,\mathbb{R}).$  If $u\in
\Lambda_{\om,K}$ and $[\om_1]$ is close enough to it, then there
exists a presymplectic $\om_2\in u$ such that $\om_1=\om_2+\eta,$
where $\eta$ is small (say, in a norm defined by a Riemannian
metric)  and vanishes near $K.$ Then $\om_2+t\eta$ is a
presymplectic path connecting $\om_2$ with $\om_1,$ constant in a
neighborhood of $K.$ This path can be used to pass from a path
connecting $\om$ with $\om_1$ to a path connecting $\om$ with
$\om_2$ and backwards. Thus if $u'$ is close enough to an element of
$\Lambda_{\om,K},$ then it also belongs to $\Lambda_{\om,K}.$

Essentially the same argument shows that $\Lambda_{\om,K}$ is
closed, since if we have an element $\om_2\in\Lambda_{\om,K}$ close
enough to $\om_1,$ then again we can use the linear segment to draw
the required path from $\om_2$ to $\om_1.$ \hfill $\square$

Theorem \ref{prehomotopy} has the following  corollary.

\begin{cor} Let $(M,\om)$ be a presymplectic manifold. Then in every class
$ a \in H^2(M,\mathbb R)$ there is a presymplectic form homotopic to
$\om$ through presymplectic forms.
\end{cor}

More difficult questions of that type are those with  Reeb foliation
involved. For example, one might ask whether for every non-zero
vector field $X$ on $M^{2n-1}$ there is a presymplectic form for
which $X$ is the Reeb field. The answer is negative even in
dimension three and the following yields a counter-example on $T^3.$

On $T^3=T^2\x S^1$ consider coordinates $(x,y,t), \ t\in [0,2].$
Take $X = T^2 \times [0,1] \subset T^2 \times S^1 \equiv T^3$ with a
coordinate system $((x,y),t)$ and define a vector field $R$ on $T^3$
by the formula $R=\cos(\pi t) \frac{\partial}{\partial t} + \sin(\pi
t) \frac{\partial}{\partial x}.$ Note that $R$ points inwards along
both boundary tori $T_0 = T^2 \times \{0\}$ and $T_1 = T^2 \times
\{1\}$ of $X.$ Assume that $R$ is a Reeb field  of a presymplectic
form $\om$ . Counting orientations of $T_0,T_1$ given by $\om$ and
$R,$ we come to contradiction with equality $\int_{T_0}\om =
\int_{T_1}\om .$

\section{Moser type theorems and contact structures}

Among classical results  of symplectic topology a special role
plays the theorem of Moser. It says that a path $\{\omega_t\}_{t \in
[0;1]}$ of cohomologous symplectic forms on a closed manifold is
an  {\it isotopy}, i.e. there exists an isotopy $\phi_t$ starting
from $\phi_0=id$ such that $\om_t=\phi_t^*\om_0.$ Let us denote an
obvious obstruction to the isotopy in presymplectic case: the Reeb
foliations should be conjugated by diffeomorphisms (see example after
Theorem \ref{mos} below). The crux of the proof is that when we write the
equation

\begin{eqnarray}\label{moser1}\sigma_t + \iota(X_t)\omega_t = 0,\end{eqnarray}

\noindent where

\begin{eqnarray}\label{moser2}\frac{d}{dt}\omega_t = d\sigma_t \end{eqnarray}

\noindent and $X_t$ generates the isotopy $\{\phi_t\}_{t \in [0;1]}$
we are looking for, then the equation has always a solution. This is
due to non-degeneration of $\om_t.$ This obviously does not extend
to the case of presymplectic forms.  However, it does extend if  we
know that presymplectic form $\sigma_t$ in \ref{moser2} satisfies
$\sigma_t(R_t)=0.$   This leads to the notion of {\it basic
cohomology} (see for example \cite{Mo}).

For any foliation $\mathfrak{F}$ consider the space $T\mathfrak{F}$
of vectors tangent to the foliation $\mathfrak{F}$ and define the
space of {\it basic n-forms} as
$$\Omega^n_b = \{\alpha \in \Omega^n \ \mid \
\forall X\in T\mathfrak{F}  \ \iota_X\alpha = \iota_X d\alpha =
0 \}.$$

The usual exterior derivative $d$ defines a mapping
$d^n_b:\Omega^n_b \rar \Omega^{n+1}_b.$ Homology groups of the
resulting chain complex are called {\it basic cohomology groups} of
$\mathfrak{F}$ and denoted $H_b^*(M;\mathfrak{F})$ (we suppress
$\mathfrak{F}$ from the notation if it is clear what foliation is
considered). Now, with a stronger assumption that the basic
cohomology class is constant, we can repeat Moser's argument to get
the following.

\begin{thm}\label{mos} Assume that $\{\om_t\}_{t \in [0,1]}$ is a path of
presymplectic forms in a fixed cohomology class on $M.$ If the Reeb
foliation is independent of $t$ and the basic cohomology class
$[\om_t]_b \in H^2_b(M)$ with respect to the Reeb foliation is
constant, then $\{\om_t\}$ is an isotopy.
\end{thm} \hfill $\square$

The assumption on basic cohomology class cannot be relaxed, since
the foliation can vary along with the parameter. For example, on
$T^3$ define $\om = dx \wedge dy, \ \alpha = sin(y) dz.$ Then for
$t$ small enough  $\om_t = \om +t d\alpha$ is presymplectic. But for
many values of $t$ the Reeb foliation of $\om_t$ is irrational
(orbits are dense), whereas for $\om_0$ it is rational (all orbits
are closed circles).

However, it is rather easy to strengthen slightly the last theorem
as follows. If two foliations are conjugated by a diffeomorphism,
then we can identify basic cohomology groups and using this
identification we can compare basic cohomology classes.

\begin{theorem} If $\om_t$ is a path of presymplectic forms on $M,$
$\phi_t$ is an isotopy  such that $\mathfrak F_t$ is conjugated to
${\mathfrak F}_0$ by  $\phi_t$ and the basic cohomology class of
$\om_t$ is constant (when we identify $H_b^2(M,\mathfrak F_t)$ with
$H_b^2(M,\mathfrak F_0)$ using $\phi_t$), then $\{\om_t\}$ is an
isotopy.
\end{theorem}\hfill $\square$

\section{Presymplectic forms of contact type}\label{contacttype}

A challenging problem concerning contact forms is to give necessary
and sufficient condition to assure that a non-degenerate 2-form on a
closed manifold can be deformed to a form of contact type. Recall
that we call a 2-form $\om$ a form of contact type if $\om = d\la ,$
where $\la$ contact. In particular, any results solve the question
of existence of contact structures. The problem was addressed by
Eliashberg  in dimension 3, and for simply connected 5-manifolds by
Geiges \cite{G2}. They both proved that there are no obstructions.

Such problems seem to be very difficult in general. We give here a
necessary condition  for presymplectic form to be of contact type.

\begin{thm} Assume that $\lambda$ is a contact form on a closed
manifold  $M$ and $\mathfrak R$ its Reeb foliation. Then
$[d\lambda]$ is non-zero in $H^2_b(M,\mathfrak R).$
\end{thm}

{\bf Proof.} Let $R$ denote the Reeb vector field of $\lambda ,$ so
that $\lambda (R) \equiv  1.$ If $[d\lambda ]_b=0,$ then there
exists $\alpha$ such that $\alpha (R)=0$ and $d\alpha = d\lambda .$
Thus  $\phi_0 = \alpha - \lambda$ is closed and equal to 1 on $R.$
It yields existence of a closed form $\phi$ which is $C^1$-close to
$\phi_0$ and such that $\ker \phi$ is integrable with compact
leaves, cf. \cite{Ti}. Since $d\lambda$ restricts to a symplectic
form on any leaf,  its cohomology class would be non-zero, which is
of course false.  \hfill$\square$

One can notice that for open manifolds such existence questions as
above  are solved positively by Gromov's h-principle (see
\cite{EM}). In particular, one can always deform a presymplectic
form to a form which is contact outside any given non-empty open set.

\section{Some examples}\label{ex}

Let $\om$ be a presymplectic form, $\om = d\al ,$ and $R$ be a
positively oriented Reeb vector field for $\om .$ Define
$f_{\al}=\al (R).$ Any other choice of $R$ corresponds to
multiplication of $f_{\al}$ by a positive function. Then we have the
following.

\begin{enumerate}
\item $\al$ is contact \ifff $f_{\al}>0$ for any choice of $R$ (and
then we can choose $R$ such that $f_{\al}$ is constant and
positive);
\item $\al$ is $R-$invariant \ifff $f_{\al}$ is constant.
\end{enumerate}

\begin{ex}
$T^2$--bundle over the circle.
\end{ex}
 Consider the following example constructed  by Carri\`{e}re in \cite{YC}.
 The manifold $T^3_A$
he considers is the $T^2-$bundle over the circle whose monodromy is
given by matrix $A \in SL(2,\mathbb Z)$ such that $tr A>2,$
$$T_A^3 = T^2 \times \mathbb R \slash (x,t) \sim (Ax,t+1).$$

\noindent Then both eigenvalues $\la, \frac{1}{\la},\ \la >1>
\frac{1}{\la}$ of $A$ are real and irrational. The eigenvectors
define two vector fields $\mu_1, \mu_2$ on $T^3_A.$ Using them we
can define two 1--forms $v_1,v_2$ on $T^2$ by the formulas
$v_1(\mu_1) = 1, \ v_1(\mu_2) = 0$ and $v_2(\mu_1) = 0, \ v_2(\mu_2)
= 1.$ This definition extends  to $T^2 \times \mathbb{R}$ by setting
$$\al_1= \la^t v_1$$
$$\al_2= \frac{1}{\la^t}v_2$$ and then gives two well
defined forms $\al_1,\al_2$ on $T_A^3.$ By direct calculation,
$d\al_1 = \ln(\la) dt \wedge\al_1$ and $d\al_2 = -\ln(\la)
dt\wedge\al_2.$ Thus  $d\al_1$ is a presymplectic form on $T^3_A$
with associated Reeb field $R = \mu_2.$ In this case, $f_{\al} = 0.$
Furthermore,
$$\phi_{\epsilon} = \al_1 + \epsilon \al_2$$  is a contact form for $\epsilon$ small
enough since
$$\phi_{\epsilon} \wedge d\phi_{\epsilon} = 2\epsilon \ln(\la) \al_1 \al_2 dt > 0.$$
This contact form is $C^1$-close to the presymplectic form $\al_1.$

Carri\`{e}re proved that the second basic cohomology group of the
Reeb bundle of $d\al_1$ vanishes, so by Theorem 5.1 there is no
contact form with the Reeb field equal to $\mu_2.$
This follows also from Taubes theorem about existence of closed
orbits of a contact Reeb field (the Weinstein conjecture).

\begin{ex}
$S^4 \times S^1$
\end{ex}

By Eliashberg's beautiful result \cite{E}, $S^1\times S^4$ has a
contact form. Namely,  if a compact almost complex manifold with
boundary of even dimension $2n,\ n>2,$ has a Morse function
(constant and maximal on the boundary) without critical points of
index greater than $n,$ then the boundary is contact. The fact that
$M$ is presymplectic is much simpler, since this is enough to notice
that $M$ is parallelizable.

We shall build directly a presymplectic form  $\om =d\al$ on $M=S^4
\times S^1$ such that the function $f_{\al}$  has both positive and
negative values, hence $\om$ cannot be $C^1$-approximated by a
contact form. Let $x_1,...,x_6$ be standard coordinates in $\mathbb
R^6.$  Consider $S^4$ as the hypersurface $(x_2-2)^2+...+x_6^2=1$ in
$\mathbb R^5 = \{ x_1=0  \}$ and rotate it around $\{ x_1=x_2=0\}.$

If we set $\al
=\frac{1}{2}(x_1dx_2-x_2dx_1+x_3dx_4-x_4dx_3+x_5dx_6-x_6dx_5),$ then
$\om = d\al = dx_1dx_2+dx_3dx_4+dx_5dx_6$ is the standard symplectic
form on $\mathbb{R}^6.$ Since $\om$ restricted to any 5-dimensional
subspace of $\mathbb{R}^6$ is presymplectic, hence $d\al=\om$
restricts to a presymplectic form on $M.$ If $X$ is orthogonal to
$M,$ then $J_{st}X,$ where $J_{st}$ is the standard almost complex
structure on $\mathbb{R}^6,$ is a Reeb field $R$ of $\om |M.$ Having
established that, it is easy to calculate that $f_{\al}$ takes both
positive and negative values.

\section{Thurston's construction}\label{thurston}

Let $p:M\rar B$ be a fibration with base $B$ and a symplectic fibre
$(F,\om_0).$ We say that the fibration is {\it symplectic} if its
structure group is $Symp(F,\om_0),$ and {\it compact symplectic} if
it is symplectic and both $B$ and $F$ are compact.  For a discussion
of symplectic fibrations see \cite{MS}, Ch.6. On each fibre
$F_b=p^{-1}(b)$ of a symplectic fibration there is  a well defined
symplectic form $\om_b$ given by $\om_0$ and local trivializations.
Denote by $i_b:F_b\subset M$ the inclusion of the fibre. The
following theorem is given by Thurston's construction, originally
applied to symplectic fibrations over a symplectic base to get
symplectic structures on total spaces \cite{Th}.

\begin{theorem}\label{thur} Let $p:M\rar B$   be a compact symplectic fibration
with fibre $(F,\om_0)$ and a presymplectic base $(B,\om_B).$ For any
class $u\in H^2(M,\mathbb R)$ such that $i_b^*u = [\om_b]$ there
exists $K>0$ and a presymplectic form $\om$ on $M$ such that
$i_b^*\om = \om_b$ and $[\om ]=u+Kp^*[\om_B].$
\end{theorem}

{\bf Proof.} Let $\{ U_{\al} , f_{\al} : p^{-1}U_{\al} \rar U_{\al}
\times F\}$ be charts of an atlas of local trivializations and
$\chi_{\al}$ a  smooth partition of unity subordinated to $\{
U_{\al}\} .$ On each set $p^{-1}U_{\al}$ we have a form $\om_{\al}$
which is pull-back of $\om$ by projection to the fibre $F.$ If
$\tau$ represents $u,$ then $\om_{\al} - \tau = d\phi_{\al}$ for
some 1-forms $\phi_{\al}.$ The formula  $\tau_1 = \tau + d\sum
(\chi_{\al}\circ p) \phi_{\al}$ defines a closed form on $M$ which
restricts to $\om_b$ on a fibre $F_b$ and represents $u.$ Then for K
large enough, $\tau_1 + Kp^*\om_B$ is a presymplectic form on $M$
representing $u+Kp^*\om_B.$ Compare \cite{Th,MS}.
\hfill$\square$

We will prove a relative version of Thurston's construction. Assume
additionally  that $j_0:N_0\subset F$ is a smooth compact
codimension zero submanifold of $F$ such  that $N_0$  is preserved
pointwise by the structure group of the given fibration. In other
words, the structure group is the group $Symp(F,N_0,\om_0)$ of
symplectomorphisms which are equal to the identity on $N_0.$ Then
$N_0$ defines submanifold $N_b$ in each $F_b$ and a submanifold
$N\subset M$ such that $N\cap F_b=N_b,\ N=N_0\times B.$ Via
pull-back of $\om_0$ under projection on $F$ one has on $N$ a form
$\om_N$ equal to $\om_0$ on fibers and zero in $B$-direction.

\begin{theorem}\label{thurl} Let $N_0\subset F$ be a compact
codimension zero submanifold and let $p:M\rar B$ be a compact
symplectic fibration with fibre $(F,\om_0)$ and structure group
$Symp(F,N_0,\om_0).$

For any presymplectic form $\omega_B$ on $B$ and any
class $u\in H^2(M,\mathbb R)$ satisfying
\begin{enumerate}
\item  $i_b^*u = [\om_b]$,

\item \label{onN}$u|N = [\om_N]$
\end{enumerate}

\noindent there exist a constant $K>0$ and  a presymplectic form $\tau$ on $M$
such that $[\tau ]=u+Kp^*[\om_B]$ and $\tau = \om_N + Kp^*\om_B$ on
$N.$
\end{theorem}

{\bf Proof}. By assumption \ref{onN}, there exist a form $\tau$
representing $u$ such that on $N$ we have  $u=\om_N.$ Let
$\phi_{\al}, \tau_1$ be forms defined in the proof of Theorem
\ref{thur}. Then $d\phi_{\al}=0$ on $N.$ Thus the form $\eta
=d\sum_{\al} (\chi_{\al}\circ p) \phi_{\al}$ on $N$ is horizontal,
since both $\eta$ and $d\eta (=0)$ vanish on vertical vectors. Thus
it is equal to $p^*\eta'$ for a closed form $\eta'$ on $B.$ Then
$\tau_1 - p^*\eta + Kp^*\om_B$ has the required properties.
\hfill$\square$

\section{Submanifolds}

Let $(M^{2n+1},\om)$ be a presymplectic manifold and let $Q^{2k+1}$ be a
smooth closed submanifold. One is tempted to define a presymplectic
submanifold as  those $Q$ that the restriction of $\om$ is
presymplectic. However, the preceding discussion shows that the Reeb
foliation plays an important role, thus it is better to impose a
condition to have it under control.

\begin{defin}
Let $(M,\om )$ be a presymplectic manifold. A submanifold $i:Q
\hookrightarrow M$ is a {\it presymplectic submanifold} if $\om_Q =
i^*\om$ is presymplectic and the Reeb bundle of $\om_Q$ is equal to
the restriction of the Reeb bundle of $\om$ to $Q$ (equivalently,
the Reeb bundle of $\om$ is tangent to $Q).$
\end{defin}

\begin{ex} If  a presymplectic form is invariant under an action of
$S^1,$ then any component of the fixed point set is a presymplectic
submanifold.
\end{ex}

For a symplectic submanifold $Q$ of a symplectic manifold $M$ and in
many other cases, the structure of a tubular neighborhood is
determined up to isotopy by $Q$ and the symplectic normal bundle.
This is not the case for presymplectic forms, since the Reeb field
might not be invariant under contractions of the neighborhood to
$Q.$ For contact submanifolds it is possible to obtain such
invariance when we allow the contact form to be multiplied by a
function. Note that this operation changes, in general, the Reeb
bundle.

Let $(M,\om )$ be a presymplectic manifold and $\mathcal R$ be the
Reeb bundle. Then the associated reduction of $TM$ to a complex
bundle is given by a choice of a bundle endomorphism $J:TM\rar TM$
such that:

\begin{enumerate}

\item $ker\, J=\mathcal R;$
\item\label{jprop2} $J^2=-Id$ on $Im\, J;$
\item\label{jprop3} $\om (JU,JV)=\om (U,V)$ for all $U\in Im\, J$ and arbitrary $V;$
\item\label{jprop4} $\om (U,JU)>0$ for all non-zero $U\in Im\, J.$
\end{enumerate}

Denote $\mathcal S=Im\, J.$ We have $TM=S\mathcal \oplus \mathcal R$
and the formula $g(U,U')=\om (U,JU')$ defines  a
Riemannian metric on $\mathcal S$. When a Reeb vector field $R$ is chosen, then we
get a Riemannian metric on $M$ defined on $\mathcal R$ by
$g(R,R)=1.$ Note that \ref{jprop3} implies that $\mathcal
S\bot\mathcal R.$ The choice of $J$ consists in a choice of a linear
complement $\mathcal S$ of $\mathcal R$ and a choice of complex
structure on it compatible with $\om$ (note that $\om$ is a linear
symplectic form on $\mathcal S).$ This gives $J$ with  required
properties if we extend by zero on $\mathcal R.$ For the resulting
Riemannian metric we have $g(JU,JV)=g(U,V)$ for all $U\in \mathcal
S$ and every $V.$

If $Q$ is a presymplectic submanifold of $(M,\om ),$ then we can
assume that $J$ preserves $TQ.$  Consider a linear complement
$\mathcal N$ of $TQ$ in $TM$ such that $\mathcal N\subset \{ U\in
TM|Q:\forall V\in TQ\ \om (U,V)=0\} .$ Since $\om$ is
non-degenerate on $\mathcal N,$  we can choose $J$ in $\mathcal N$
compatible with $\om .$ Then, with respect to the Riemannian metric
$g$ (on the whole M) constructed as above, $\mathcal N$ is an
orthogonal complement to $TQ$ and $g$ defines a horizontal
distribution $\mathcal H\subset T(\mathcal N).$ We get now a
presymplectic form on the total space of $\mathcal N,$ linear in
every fibre, by the formula

$$ \om^{\mathcal N}_W(V_1+H_1,V_2+H_2)=
\om_{p(W)}(p_*H_1,p_*H_2)+\om_{p(W)}(V_1,V_2),$$

\noindent where $W\in \mathcal N,\ V_1,V_2$ are vertical, $H_1,H_2$
horizontal vectors in $T_W(\mathcal N)$ and $p:\mathcal N\rar Q$ is
the bundle projection of $TM$ restricted to $\mathcal N.$

In this way we get complex structures on the normal bundle $\nu
Q=(TM|Q){/}TQ$ and compatible with the presymplectic forms on the
total space of $\nu Q$ via isomorphisms $\mathcal N\subset TM|Q\rar
\nu Q.$ Any two such complex structures obtained in this way are
isomorphic.

Consider now the exponential map $exp: \nu Q\rar M$ given by the
metric $g.$ It is a diffeomorphism near $Q$, thus $\om^{\mathcal N}$
defines a non-degenerate form linearized in the normal direction
which we also denote $\om^{\mathcal N},$ in a neighborhood of $Q$.
The form restricted to $TM|Q$ is equal to $\om ,$ thus these
forms in a small neighborhood of $Q$ are close enough to have a linear
segment contained in the space of non-degenerate
forms which connects them. Hence, for some open tubular neighborhoods $Q\subset U\subset
U_1$ such that $\overline U_1-U\cong \partial \overline U\x [0,1],$
there is a smooth non-degenerate form $\om'$ such that:

\begin{enumerate}

\item $\om'=\om$ on $M-U_1;$
\item $\om'=\om^{\mathcal N}$ on $U;$
\item $\om'|\overline U_1-U$ is a smooth linear combination  $\la(t)\om
|\partial\overline U_1 + (1-\la(t))\om^{\mathcal
N}|\partial\overline U,$ where $ t\in [0,1]$ and $\la$ is an
appropriate smooth function changing from $0$ to $1.$

\end{enumerate}

With our choices, $t\om +(1-t)\om'$ is a path of non-degenerate
forms. It connects $\om$ with $\om'$ and is equal to $\om$ outside
$U_1.$

Let us assume that we have two presymplectic forms $\om_0, \om_1$
which coincide on $Q$ and define isomorphic structures on the normal
bundle $\nu Q.$ Then for any choice of normal subbundles $\mathcal
N_0, \mathcal N_1$ and horizontal distributions $\mathcal H_0,
\mathcal H_1$ there is a path $(\mathcal N_t, \mathcal H_t, J_t);\
t\in [0,1]$ connecting $(\mathcal N_0,  \mathcal H_0, J_0)$ with
$(\mathcal N_1, \mathcal H_1,J_0).$ Thus the forms $\om'_0, \om'_1$
are homotopic through a path $\om'_t$ of non-degenerate forms.

 This in turn gives a path of non-degenerate forms
connecting $\om_0$ with a presymplectic form equal to $\om_0$
outside $U_1$ and to $\om_1$ in a neighborhood of $Q.$ It is defined
as follows. Take tubular neighborhoods $Q\subset U'\subset
U_2\subset U$ such that $U_2-U'$ is a product by an interval.
Construct, using $\om'_t,$ a path of non-degenerate forms equal to
$\om_0$ outside $U_1,$ with the end form $\eta'_1$ equal to $\om'_1$
in $U'.$ Next deform $\eta'_1$ to a form $\eta_1$ equal to $\om_1$
in $U',$ with the help of the path connecting $\om'_1$ with $\om_1.$
Finally, there  is  a presymplectic homotopy connecting $\om_0$ with
a presymplectic form equal to $\om_0$ on $M-U'_1$ and to $\om_1$
on $U'',$ where $U'_1,U''$ are tubular neighborhoods of $Q$ such
that $U''\subset U', \ U_1\subset U'_1.$ This can be done in two
steps. First deform $\eta_1$ to a presymplectic form with the
required properties by Theorem \ref{dusa}, then from obtained
path of non-degenerate forms pass to a presymplectic homotopy using
Theorem \ref{prehomotopy}, relative to $K=M-U'_1.$

This proves the following tubular neighborhood theorem.

\begin{thm}\label{deform} Assume that $\om_0, \om_1$ are
presymplectic forms such that $Q$ is a presymplectic submanifold
with respect to both of them. If the presymplectic forms coincide on
$Q$ and the complex normal vector bundles to $Q$ induced by $\om_1$
and $\om_2$ are isomorphic, then there exist tubular neighborhoods
$U_0, U_1$ of $Q,\ U_0\subset U_1,$  and  a presymplectic homotopy
connecting $\om_0$ to a presymplectic form equal to $\om_0$  on
$M-U_1$ and to $\om_1$ on $U_0.$
\end{thm} \hfill $\square$

\section{Starlike structures}

An open book decomposition of $M$ is defined by a quadruple $\{B,P,\pi
,\phi\}$ where $B$ is a codimension 2 submanifold of $M,$ $\pi
:M-B\rar S^1$ is a smooth fibration with fibre $Int\, P$ and $\phi$ is the
gluing diffeomorphism of the fibration. It is assumed that $\partial
P = B$ and $\phi$ is equal to the identity on $U-B,$ where $U$ is
a neighborhood  of
$\partial P,$ so that the normal bundle of  $B$ is trivial. If $U$
is a collar, then a neighborhood of $B$ is diffeomorphic to $B\x D^2.$
$P$ is called {\it page} of the decomposition and $B$ it's {\it binding.}

Consider an open book decomposition of $M$ satisfying the following
additional assumptions.

\begin{enumerate}
\item $P$ is endowed with a  symplectic form $\om_0,$
\item $\phi$ preserves $\om_0.$
\end{enumerate}

Then the fibration $M-B\rar S^1$ is symplectic and one can apply
Thurston type construction described in Theorem \ref{thurl}
to prove the following.

\begin{theorem}
Consider a closed smooth manifold $M^{2n+1}$ and an open book
decomposition of $M$ with symplectic page $(P,\om_0 )$ and the
gluing diffeomorphism  $\phi \in Symp(P,U,\om_0).$ For any
cohomology class $u\in H^2(M,\mathbb R)$ such that $i^*_tu=[\om_t] $
there exists a presymplectic form $\om$ on $M$ such that:

\begin{enumerate}
\item $[\omega]=u;$
\item outside of a neighborhood of the binding $B$ the Reeb
vector field is transversal to pages;

\item the binding $B$ is a presymplectic submanifold of $M.$

\end{enumerate}
\end{theorem}

We shall prove this for a more general structure of
 {\it starlike decomposition} of $M.$ By that we mean  a
 quadruple $\{ C,S,\pi ,\phi \} ,$ where $C$ is a
 codimension $2k$ submanifold of $M,$ $\pi
:M-C\rar S^{2k-1}$ is a smooth fibration with fibre $Int\, S$ and
the gluing map $\phi : S^{2k-2}\rar Diff(S,U)$ of the fibration
takes values in the group of diffeomorphisms equal to the identity
on an open neighborhood $U$ of $\partial S.$ It is assumed also that
$\partial S = C$ and thus the normal bundle of $C$ is trivial. We
will call  $S$  {\it spine} of the decomposition and $C$ it's {\it
core.}

Assume now  that $S$ is endowed with a  symplectic form $\om_0$ and the
gluing map has values in the symplectomorphism group $Symp(S,U,\om_0),$ so
that a neighborhood of $C$ is diffeomorphic to $C\x Int\, D^{2k}.$
In the sequel we assume that the neighborhood $U$ in the definition
 is chosen in that  way.
As above, $S_t=\pi^{-1}(t),$  $i_t:S_t\subset M$
is the inclusion, $\omega_t$ is the  symplectic form induced
by $\om_0$ on $S_t, t\in S^{2k-1}.$

\begin{theorem}
Consider a closed smooth manifold $M^{2n+1}$ and a starlike
decomposition of $M$ with a symplectic spine $(S,\om_0 )$ and a
gluing map  $\phi : M-C\rar  Symp(S,U,\om ).$ For any
cohomology class $u\in H^2(M,\mathbb R)$ such that $i_t^*u=[\om_t] $ and
an open neighborhood $V$ of $C,$
there exists a presymplectic form $\om$ on $M$ such that:

\begin{enumerate}
\item $[\omega ]=u;$
\item  outside of a tubular neighborhood $W\subset V$ of the core $C$ the Reeb
vector field is transversal to spines;

\item the core $C$ is a presymplectic submanifold of $M.$

\end{enumerate}
\end{theorem}

{\bf Proof.} Let us start with the following observation. Consider
a manifold $X$ with boundary. If $\eta$ is a symplectic or
non-degenerate  form on $X,$ then one can deform $\eta$ near the
boundary to get a non-degenerate (but perhaps not closed) form which
is of the form $p^*\eta_0 + dt\wedge p^*\phi$ on a collar
$W=\partial X\x [0,1)$ of the boundary.  Here $p:W\rar\partial X$ is
the projection and $\eta , \phi$ are some forms on $\partial X.$  We
will use this to construct a non-degenerate form on $M$ and to
define a presymplectic form by Theorem \ref{dusa}.

Let a tubular neighborhood $W_1\subset U$ of the core
diffeomorphic to $C\x Int\, D^{2k}$ be given. Applying the remark
above to the
standard symplectic form on $\mathbb D^{2k}$ , we get on $D^{2k}$  a
non-degenerate 2-form $\alpha$ equal to the standard symplectic form
in a neighborhood of $0$ and product near the boundary, i.e. equal
to $dt\wedge\lambda +d\lambda$ on $(1-\epsilon ,1+\epsilon)\x S^{2k-1},$
where $\lambda$ is a standard contact structure on $S^{2k-1}.$  Here
$\{ 1+\epsilon\} \x S^{2k-1}$ is the boundary,
$t\in (1-\epsilon ,1+\epsilon].$ Thus on $C\x D^{2k}$ we
have a presymplectic form $\Omega_1=\om_0|C\x \alpha$ which has
required properties near $C.$

Let $W_2\subset W_1$ be a subcollar corresponding to
the interval $[0,1)\subset [0,1+\epsilon ).$
The manifold $M'=M-W_2$ is diffeomorphic to the fibre bundle over
$S^{2k-1}$ with fibre $S'=S-(S\cap W_2)$ provided by the starlike
decomposition. Let $\Omega_3$ be a presymplectic form constructed in
Theorem \ref{thurl} on $M-W_2$ such that $i_t^*\Omega_3=\om_t$ and equal
to the product $\om_0\times d\lambda$ in $W_1-W_2.$
 As above, by deforming $\om_0$ near $C=\partial S$ we can
assume that $\om_0$ is of the form $\om_0|C+dt\wedge\eta ,$ for a
form $\eta$ on $C.$ Since the resulting form is non-degenerate,
$\eta$ is non-zero on the Reeb
subbundle in $TC.$ So we have now a non-degenerate form which is
presymplectic outside a neighborhood of the boundary and product
near the boundary.

Now let $A\subset W_1$ be a subcollar corresponding to
the interval $[0,1-\epsilon')\subset [0,1+\epsilon ),$ where
$\epsilon'<\epsilon ,$ and $B$ a subcollar corresponding to
$[0,1+\epsilon').$  Consider $\Om_1$ on $A$ and
$\Om_3$ on $M-B.$
We connect these two forms by a form on
$B-\overline A=C\x (a,b)\x S^{2k-1}$ as follows.
Define $\Omega_2=\om_0|C \x dt\wedge (h(t)\lambda
+(1-h(t))\eta ) \x d\lambda ,$ for appropriately chosen smooth
function $h:(a,b)\rar \er$ vanishing near $a$ and
equal to 1 near $b.$ Notice that
for Reeb fields $R_1=ker\, {\omega_0|C}, R_3=ker\, {d\lambda}$ we have
$\lambda (R_1)=\eta (R_3)=0,$ so that $ker \Omega_2$ is of dimension
1 spanned by $(h-1)R_3+hR_1.$  This implies that $\Omega_2$ is
non-degenerate and we get a smooth non-degenerate form on $M.$ Deforming it
according to (the relative version of) Theorem \ref{dusa} we get required
presymplectic structure.\hfill$\square$

\section{Presymplectic surgery}

The classical surgery on smooth manifolds is performed by cutting a
product neighborhood $S^k\x D^{n-k}$ of an embedded  sphere $S^k$
and gluing $D^{k+1}\x S^{n-k-1}$ along the boundary. This operation
can be described as deleting $S^k\x D^{n-k}$ from $M^n$ and
$S^n=S^k\x D^{n-k}\cup D^{k+1}\x S^{n-k-1}$ and gluing the
manifolds obtained in this way along the obvious diffeomorphism of
boundaries. In the present context, prominent examples of
constructions of that type are given by Eliashberg for contact, and
by Gompf for symplectic manifolds. Gompf assumes that two symplectic
submanifolds $V,V'$ of codimension two embedded in respectively
$M,M',$ are symplectomorphic and that the normal bundles $\nu,\nu'$
of these submanifolds satisfy $c_1(\nu)+c_1(\nu')=0.$ Then one  can
perform surgery along tubular neighborhoods of these submanifolds
resulting in a new symplectic manifold. We want to settle a
presymplectic analog of Gompf's construction.

Consider two presymplectic submanifolds $V \subset (M,\om), \ V'
\subset (M',\om').$ We assume that there exists a presymplectic
diffeomorphism $g: V \rightarrow V'$ and an orientation reversing
linear mapping $G:\nu V \rar \nu V'$ of normal bundles covering $g.$
Using some auxiliary Riemannian metrics on $M,M'$ we identify total
spaces of the normal unit disk bundles with neighborhoods $N,N'$ of
$V$ and $V'$. We can assume that $G$ is an isometry of normal
bundles. Let $N_{\delta},N_{\delta}'$ correspond to $\delta -$disk
bundles. Fix $0<\epsilon <\frac12.$ $G$ induces a diffeomorphism of
$N_{1-\eps}$ to $N_{1-\eps}'$ which we also denote by $G.$ Glue
$M-N_{\eps}$ with $M'-N_{\eps}$ along boundaries using $G.$ We claim
that the resulting manifold, denoted by $M\cup_{g}M'$ admits  a
presymplectic structure.

\begin{theorem}\label{surgery} Under the above assumptions, there
exists a presymplectic form on $M\cup_{g}M',$  equal to the given
presymplectic forms in $M-N \cup M'-N' \subset M\cup_{g}M'.$
\end{theorem}

{\bf Proof.} As in Section \ref{symplhom}, take a 1-form $\eta$ on
$M$ such that $\eta(R)=1,$ where $R$ is a Reeb field of $\om ,$ and
define a symplectic form $\phi = p^*\om + \delta (d\theta\wedge\eta
- \theta\, d\eta)$ on $M\x \mathbb R.$ When we restrict $\phi$ to
$M\x\{ 1\},$ we get a presymplectic form $\om + d\eta $ homotopic to
$\om ,$ if $\delta$ is small enough. By Theorem \ref{deform}, there
is a deformation of $\om',$ supported in $N'_{1-\eps},$ to a
presymplectic form $\om''$ equal to $(G^{-1})^*(\om +d\eta )$ in
$N'_{\eps}.$

Consider the disjoint sum $$Z = M \sqcup (U \times \mathbb{R})
\sqcup M',$$ and identify

$$M' \supset N'_{1-\eps} \backepsilon x \rar (G^{-1}(x),1)
\in U \times \{1\} \subset U \times \mathbb{R}.$$

In $N_{1-\eps}\x \mathbb R$ consider a  hypersurface $(M - N_{\eps})
\cup (\partial N_{1-\eps} \times [0,1]) \cup (M' - N'_{\eps}),$ where
on the level $\theta = 1$ we use the above identification.

This is a standard exercise that one can smooth out the corners of
the hypersurface we defined, so we get a smooth hypersurface  in a
symplectic manifold, hence a presymplectic form on it. Now a
scrutiny of forms and identifications shows that we get a
presymplectic form satisfying the requirements of the theorem. \bk

Classical surgery can be used to  simplify manifolds, for instance
to construct a simply connected manifold cobordant to any given
oriented one. As an application of Theorem \ref{surgery} we will
show that one can get a 1-connected manifold out of a given
presymplectic manifold by presymplectic surgeries on a number of
circles.

\begin{prop} If each generator of the fundamental group of a
presymplectic manifold $M$ can be represented by a closed orbit of
the Reeb foliation, then presymplectic surgeries on these circles
transform $M$ into a simply connected presymplectic manifold.
\end{prop}

{\bf Proof.} Take $S^{2n+1}=S^1\times
D^{2n}\cup_{\text{id}_{\partial}}D^2 \times S^{2n-1}.$ It admits a
presymplectic form such that $S^1\times \{ 0\}$ is a closed orbit of
the Reeb foliation. Since presymplectic manifolds are orientable by
definition, thus normal bundles of any embedded  circle
is trivial.Therefore for each closed orbit of the Reeb field one can
perform the surgery of $M$ with $(M',V')=(S^{2n+1},S^1\times \{
0\}).$ Notice that we perform in fact the classical surgery on a
1-sphere so that the homotopy class of the circle is killed. The
proposition follows.\bk

In order to be able to apply the last proposition we  need a simple
lemma. Denote the space of linear forms of maximal rank on $\mathbb
R^k$ by $\Om(\mathbb{R}^{k}).$ This means that
$\Om(\mathbb{R}^{2n})$ is the space of symplectic linear forms and
$\Om(\mathbb{R}^{2n+1})$ is the space of presymplectic linear forms.
Moreover, let $\Om^+(\mathbb{R}^{2n})$  denote the component of
symplectic forms compatible with the orientation of $\mathbb
R^{2n}.$ Then the following holds.

\begin{lemma}\label{simp_conn_om} The space
$\Om(\mathbb{R}^{2n+1})$ is simply connected.
\end{lemma}

{\bf Proof.} There exists a  fibration $\Om^+(\mathbb R^{2n})
\mapto{i} \Om(\mathbb{R}^{2n+1}) \mapto{\pi} S^{2n},$ defined by
$\pi (\om ) = R,$ where $R$ is the unit Reeb vector of  $\om$
compatible with the orientation of $\mathbb R^{2n+1}$ and the
orientation  provided by $\om $ on the orthogonal complement of the
Reeb line. Since $\Om^+(R^{2n})$ is simply connected (\cite{MS}, Ch.
2), hence so is $\Om(\mathbb{R}^{2n+1}).$ \hfill $\square$

If $S^1 \hookrightarrow M$ is any embedded circle and $S^1 \times
D^{2n}$ is a tubular neighborhood of this circle, then by Lemma
\ref{simp_conn_om} we can assume that on a smaller tubular
neighborhood $S^1 \times D_1^{2n}$ our form $\om$ is the pullback
$\pi^*\om_{st}$ (where $\om_{st}$ denotes the standard symplectic
form on $D_1^{2n} \subset \mathbb{R}^{2n}$). Finally, Theorem
\ref{dusa} shows that one can always find a presymplectic form
which enables to apply presymplectic surgery to kill the fundamental
group.

\begin{prop}\label{closedorbit} If $M$ is a
presymplectic manifold, then for every cohomology class $a \in
H^2(M,\mathbb{R}),$ any given homotopy class of non-degenerate forms
and arbitrary elements $g_1,\ldots,g_n \in \pi_1(M),$ there exists a
pre\-symp\-lec\-tic form $\om$ in the given homotopy class such that
$[\om]=a$  and there exists embedded circles representing
$g_1,\ldots,g_n$ which are closed orbits of the Reeb foliation.
\hfill$\square$
\end{prop}

\begin{rem} {\em By Theorem \ref{prehomotopy}, we have also the following result.
If we are given a presymplectic form on $M$ (under assumptions of
Proposition \ref{closedorbit}) in a prescribed homotopy class, we
can deform it, through presymplectic forms, to a form satisfying all
the requirements of Proposition \ref{closedorbit}.}
\end{rem}

\bibliographystyle{amsalpha}

\medskip

\noindent {\bf Mathematical Institute, Wroc\l aw University,

\noindent pl. Grunwaldzki 2/4,

\noindent 50-384 Wroc\l aw, Poland}

\medskip

\NI and

\medskip

\noindent {\bf Department of Mathematics and Information Technology,

\noindent University of Warmia and Mazury,

\noindent \.{Z}o\l nierska 14A, 10-561 Olsztyn, Poland}

\medskip

\begin{flushleft}
\tt hajduk@math.uni.wroc.pl
\end{flushleft}
\medskip

\noindent {\bf West Pomeranian Technological University,

\noindent Mathematical Institute}

\noindent {\bf Al. Piast\'{o}w 48/49, 70--311 Szczecin, Poland}
\medskip

\begin{flushleft}
 \tt rafal\_walczak2@wp.pl
\end{flushleft}

\end{document}